\documentclass{amsart}
\usepackage{amsmath}
\usepackage{amsfonts}
\usepackage{amssymb}
\usepackage{verbatim}
\usepackage{calc}
\usepackage{amscd}

\usepackage{color}
\definecolor{darkblue}{rgb}{0, 0, .4}
\definecolor{grey}{rgb}{.7, .7, .7}
\usepackage[breaklinks]{hyperref}
\hypersetup{
	colorlinks=true,
	linkcolor=darkblue,
	anchorcolor=darkblue,
	citecolor=darkblue,
	pagecolor=darkblue,
	urlcolor=darkblue,
	pdftitle={},
	pdfauthor={}
}

\usepackage{graphicx}
\usepackage{graphics}
\usepackage{epsfig}

\newtheorem{theorem}{Theorem}[subsection]
\newtheorem{corollary}[theorem]{Corollary}
\newtheorem{proposition}[theorem]{Proposition}
\newtheorem{lemma}[theorem]{Lemma}
\theoremstyle{definition}
\newtheorem{definition}[theorem]{Definition}
\newtheorem{example}[theorem]{Example}

\theoremstyle{remark}
\newtheorem{remark}[theorem]{Remark}
\numberwithin{equation}{section}

\renewcommand{\a}{\alpha}
\newcommand{\e}{\varepsilon}

\author[C. Berg]{Chris Berg}
\author[B. Jones]{Brant Jones}
\author[M. Vazirani]{Monica Vazirani}

\thanks{The first author was supported in part by NSF grant DMS-0135345.  The second author was supported in part by NSF grant DMS-0636297.  The third author was supported in part by NSF grant DMS-0301320}
\address{Department of Mathematics, UC Davis, USA, 95616}
\email{[berg, brant, vazirani]@math.ucdavis.edu}

\newlength\cellsize 
\savebox2{%
\begin{picture}(18,18)
\put(0,0){\line(1,0){18}}
\put(0,0){\line(0,1){18}}
\put(18,0){\line(0,1){18}}
\put(0,18){\line(1,0){18}}
\end{picture}}
\newcommand\cellify[1]{\def\thearg{#1}\def\nothing{}%
\ifx\thearg\nothing
\vrule width0pt height\cellsize depth0pt\else
\hbox to 0pt{\usebox2\hss}\fi%
\vbox to 18\unitlength{
\vss
\hbox to 18\unitlength{\hss$#1$\hss}
\vss}}
\newcommand\tableau[1]{\vtop{\let\\=\cr
\setlength\baselineskip{-16000pt}
\setlength\lineskiplimit{16000pt}
\setlength\lineskip{0pt}
\halign{&\cellify{##}\cr#1\crcr}}}
\savebox3{%
\begin{picture}(15,15)
\put(0,0){\line(1,0){15}}
\put(0,0){\line(0,1){15}}
\put(15,0){\line(0,1){15}}
\put(0,15){\line(1,0){15}}
\end{picture}}
\newcommand\expath[1]{%
\hbox to 0pt{\usebox3\hss}%
\vbox to 15\unitlength{
\vss
\hbox to 15\unitlength{\hss$#1$\hss}
\vss}}

\begin{document}

\title[A bijection on core partitions]{A bijection on core partitions and a parabolic quotient of the affine symmetric group}

\begin{abstract}
Let $\ell,k$ be fixed positive integers.  In \cite{BV}, the first and third authors established a bijection between $\ell$-cores with first part equal to $k$ and $(\ell-1)$-cores with first part less than or equal to $k$.  This paper gives several new interpretations of that bijection.  The $\ell$-cores index minimal length coset representatives for $\widetilde{S_{\ell}} / S_{\ell}$ where $\widetilde{S_{\ell}}$ denotes the affine symmetric group and $S_{\ell}$ denotes the finite symmetric group.  In this setting, the bijection has a beautiful geometric interpretation in terms of the root lattice of type $A_{\ell-1}$.  We also show that the bijection has a natural description in terms of another correspondence due to Lapointe and Morse \cite{LM}.
\end{abstract}

\date{\today}

\maketitle

\section{Introduction}

Core partitions are combinatorial objects that appear naturally in various type $A$ settings.  They are used in the modular representation theory of the symmetric group to describe the blocks \cite{JK}.  In the geometry of the affine Grassmannian, cores index Schubert varieties and related homology classes called $k$-Schur functions \cite{LM}.  At the level of Coxeter groups, cores correspond to minimal length coset representatives for the parabolic quotient $\widetilde{S_{\ell}} / S_{\ell}$ where $\widetilde{S_{\ell}}$ denotes the affine symmetric group and $S_{\ell}$ denotes the finite symmetric group.

The impetus for this paper is combinatorial.  In \cite{BV}, the first and third authors showed that the number of $\ell$-cores with first part equal to $k$ is $\binom{k+\ell -2}{k}$ using a bijection $\Phi_{\ell}^k$ between $\ell$-cores with first part $k$ and $(\ell-1)$-cores with first part less than or equal to $k$.  The fact that such a projection exists is remarkable in part because the Coxeter group of affine type $A_{\ell-2}^{(1)}$ is not a parabolic subgroup of $A_{\ell-1}^{(1)}$.  

In this work we review some combinatorial models for cores and interpret the bijection in various guises.  In the Coxeter system setting, the bijection has a geometric interpretation as a projection from the root lattice of type $A_{\ell-1}$ to an embedded copy of the root lattice of type $A_{\ell-2}$; see Figure~\ref{f:sl3}.  We observe that the bijection reduces the Coxeter length of the corresponding minimal length coset representative by exactly $k$.  We also show that the bijection has a natural description in terms of another correspondence between $\ell$-cores and $(\ell-1)$-bounded partitions due to Lapointe and Morse \cite{LM}.

\subsection{Organization}

In Section~\ref{s:intro} we introduce the bijection $\Phi_{\ell}^k$ in terms of partition diagrams.  In Section~\ref{s:geometry}, we review the correspondence between $\ell$-cores and minimal length coset representatives for $\widetilde{S_{\ell}} / S_{\ell}$ where $\widetilde{S_{\ell}}$ denotes the affine symmetric group and $S_{\ell}$ denotes the finite symmetric group.  In Section~\ref{s:geometric_bijection}, we give a geometric version of the bijection $\Phi_{\ell}^k$ on the root lattice of type $A_{\ell-1}$.  In Section~\ref{s:l-m}, we show that $\Phi_{\ell}^k$ also has a natural description in terms of bounded partitions using the correspondence
\[ \rho_{\ell-1} : \{ \ell\text{-cores} \} \rightarrow \{ \text{partitions with first part } \leq \ell-1 \} \]
due to Lapointe and Morse \cite{LM}.

\section{Definitions, Notation and a Review of the Bijection}\label{s:intro}

\subsection{Preliminaries}
Let $\lambda = (\lambda_1, \ldots, \lambda_r)$ be a partition of $n$ and $\ell \geq 2$ be an integer.  Throughout this paper, all of our partitions are drawn in English notation.  We will use the convention $(x,y)$ to denote the box which sits in the $x^{\textrm{th}}$ row and the $y^{\textrm{th}}$ column of the Young diagram of $\lambda$.  We sometimes abuse notation and refer to row $i$ in the diagram of $\lambda$ as $\lambda_i$.  $\mathcal{P}$ will denote the set of all partitions. The length of a partition $\lambda$ is defined to be the number of nonzero parts of $\lambda$ and will be denoted $len(\lambda)$.  For a fixed choice of $\ell$, the \em residue \em of box $(i,j)$ is defined to be the least nonnegative integer $r(i,j)$ congruent to $j-i$ modulo $\ell$.  We let $\delta_{i,j}$ denote the Kronecker delta function.  The \textit{hook length} of the $(a,c)$ box of $\lambda$ is defined to be the number of boxes to the right and below the box $(a,c)$, including the box $(a,c)$ itself. It will be denoted \textit{$h_{(a,c)}^{\lambda}$}.

\begin{definition}\label{d:lcore}
A partition $\lambda$ is an \em $\ell$-core \em if for every box $(a,b)$ in the Young diagram of $\lambda$, we have $\ell \nmid h_{(a,b)}^{\lambda}$. 
\end{definition}

The set of all $\ell$-cores will be denoted $\mathcal{C}_{\ell}$. The subset of $\mathcal{C}_{\ell}$ having first part $k$ will be denoted $\mathcal{C}_{\ell}^k$ and the subset of $\mathcal{C}_{\ell}$ having first part $\leq k$ will be denoted $\mathcal{C}_{\ell}^{\leq k}$. See James and Kerber's book \cite{JK} for more background on partitions and $\ell$-cores.  Definition~\ref{d:lcore} is most useful for our purposes; $\ell$-cores are more commonly defined as partitions having no removable $\ell$-rim hook.

\subsection{$\beta$-numbers and abaci} The notion of $\beta$-numbers can be found in \cite{JK}. Here we give a modified description of the $\beta$-numbers.

Each partition $\lambda = (\lambda_1, \dots, \lambda_r)$ is determined by its hook lengths in the first column, i.e. the $h_{(i,1)}^{\lambda}$. From a sequence $(\alpha_1, \dots, \alpha_r)$ of positive decreasing integers one obtains a partition $\mu$ by requiring that the hook length $h_{(i,1)}^{\mu} = \alpha_i$ for $1 \leq i \leq r$. This gives a bijection between the set of partitions and the set of strictly decreasing sequences of positive integers. 

One can generalize this process by looking at the set $B$ of infinite sequences $b = (b_1, b_2, \dots )$ of integers. We give $B$ the group structure of component-wise addition. We define the element \textbf{1} $= (1,1,1,1, \dots) \in B$. Let $S$ denote the subgroup generated by \textbf{1} under addition, so $S = \{ (n,n,n, \dots ) : n \in \mathbf{Z} \}$. A sequence $b = (b_1, b_2, \dots ) \in B$ is said to \textit{stabilize} if there exists an $n$ so that $b_i - b_{i+1} = 1$ for all $i > n$. The set $\mathcal{B}$ is defined to be the subset of $B$ of strictly decreasing sequences that stabilize, modulo the added relation $\equiv$ that two sequences are equivalent if their difference is in $S$.
\begin{example}
$(11, 7, 4, 1, -1, -2, -3, \dots)$ is in $\mathcal{B}$. In $\mathcal{B}$, we have 
\[ (9, 5, 2, -1, -3, -4, -5, \dots) \equiv (11, 7, 4, 1, -1, -2, -3, \dots). \]
\end{example}

We define a bijection $\beta$ between the set $\mathcal{P}$ of partitions and $\mathcal{B}$. To a partition $\lambda = (\lambda_1, \lambda_2, \dots \lambda_r)$ of length $r$, we define $\beta(\lambda)$ to be the equivalence class of $(h_{(1,1)}^{\lambda}, h_{(2,1)}^{\lambda}, h_{(3,1)}^{\lambda}, \dots , h_{(r,1)}^{\lambda}, -1, -2, -3, -4, \dots)$ in $\mathcal{B}$.  

\begin{example}
$\beta(8,5,3,1)$ is the equivalence class of $(11, 7, 4, 1, -1, -2, -3, \dots).$
\end{example}

An \textit{abacus diagram} is a diagram containing $\ell$ columns labeled $0, 1, \dots, \ell-1$, called \textit{runners}.  The horizontal cross-sections or rows will be called \em levels \em and runner $i$ contains entries labeled by $r \ell + i$ on each level $r$ where $-\infty < r < \infty$.  We draw the abacus so that each runner is vertical, oriented with $-\infty$ at the top and $\infty$ at the bottom, with runner 0 in the leftmost position, increasing to runner $\ell-1$ in the rightmost position.  Entries in the abacus diagram may be circled; such circled elements are called \textit{beads}. Entries which are not circled will be called \textit{gaps}.  The linear ordering of the entries given by the labels $r \ell + i$ is called the \em reading order \em of the abacus and corresponds to scanning left to right, top to bottom.

\begin{example}
The following abacus diagram has beads in positions ($\ldots$, -3, -2, -1, 1, 2, 4, 5, 8) and gaps in positions (0, 3, 6, 7, 9, 10, 11, $\dots$). Level 0 is the row which contains $0,1,2$. 
\noindent
\\
\\

\begin{center}
\begin{picture}(80,120)(0,10)
\put (42.5,96){.}
\put (42.5,92){.}
\put (42.5,100){.}
\put (72.5,96){.}
\put (72.5,92){.}
\put (72.5,100){.}
\put (11.5,96){.}
\put (11.5,92){.}
\put (11.5,100){.}
\put (9,80) {-3}
\put (39,80){-2}
\put (69,80){-1}
\put (10,65) {0}
\put (40,65){1}
\put (70,65){2}
\put (10,50){3}
\put (40, 50){4}
\put (70,50){5}
\put (10,35){6}
\put (40,35){7}
\put (70,35){8}
\put (10,20){9}
\put (38,20){10}
\put (68,20){11}
\put(0,105){\line(1,0){85}}
\tiny
\put (0,120){Runner}
\put (30,120){Runner}
\put (60,120){Runner}
\normalsize
\put (10,110){0}
\put (40,110){1}
\put (70,110){2}
\put (-40,65){Level 0 $\to$}
\put (-42,80){Level -1 $\to$}
\put (-40,50){Level 1 $\to$}
\put (-40,35){Level 2 $\to$}

\put (72.5,53){\circle{13}}
\put (72.5,68){\circle{13}}
\put (72.5,38){\circle{13}}
\put (42.5,68){\circle{13}}
\put (42.5,53){\circle{13}}

\put (72.5,83){\circle{13}}
\put (42.5,83){\circle{13}}
\put (12.5,83){\circle{13}}
\put (42.5,12){.}
\put (42.5,8){.}
\put (42.5,4){.}
\put (72.5,12){.}
\put (72.5,8){.}
\put (72.5,4){.}
\put (11.5,12){.}
\put (11.5,8){.}
\put (11.5,4){.}
\end{picture}
\end{center}
\end{example}

A representative $\omega$ of $\beta(\lambda)$ will be called a \textit{set of $\beta$-numbers} for $\lambda$.  Suppose $\omega = (\omega_1, \omega_2, \dots )$ is a set of $\beta$ numbers for $\lambda$. An \em abacus \em for $\lambda$ is obtained by circling the entries of $\omega$ in an abacus diagram. 

\begin{example}
The following two diagrams are abaci for $\lambda = (8,5,3,1)$, the first comes from the $\beta$-numbers $(11,7,4,1,-1,-2,-3 , \dots)$ and the second comes from the equivalent $\beta$-numbers $(9,5,2,-1,-3,-4,-5, \dots )$.  We list the beads in reverse reading order to be compatible with stability in $\mathcal{B}$.
\\\\
$$
\begin{array}{lccccccr}
\begin{picture}(80,80)
\put (42.5,96){.}
\put (42.5,92){.}
\put (42.5,100){.}
\put (72.5,96){.}
\put (72.5,92){.}
\put (72.5,100){.}
\put (11.5,96){.}
\put (11.5,92){.}
\put (11.5,100){.}
\put (9,80) {-3}
\put (39,80){-2}
\put (69,80){-1}
\put (10,65) {0}
\put (40,65){1}
\put (70,65){2}
\put (10,50){3}
\put (40, 50){4}
\put (70,50){5}
\put (10,35){6}
\put (40,35){7}
\put (70,35){8}
\put (10,20){9}
\put (38,20){10}
\put (68,20){11}
\put (72.5,23){\circle{13}}
\put (42.5,68){\circle{13}}
\put (42.5,53){\circle{13}}
\put (42.5,38){\circle{13}}

\put (72.5,83){\circle{13}}
\put (42.5,83){\circle{13}}
\put (12.5,83){\circle{13}}
\put (42.5,12){.}
\put (42.5,8){.}
\put (42.5,4){.}
\put (72.5,12){.}
\put (72.5,8){.}
\put (72.5,4){.}
\put (11.5,12){.}
\put (11.5,8){.}
\put (11.5,4){.}
\end{picture}
& & & & & & &
\begin{picture}(80,80)
\put (42.5,96){.}
\put (42.5,92){.}
\put (42.5,100){.}
\put (72.5,96){.}
\put (72.5,92){.}
\put (72.5,100){.}
\put (11.5,96){.}
\put (11.5,92){.}
\put (11.5,100){.}

\put (9,80) {-6}
\put (39,80){-5}
\put (69,80){-4}
\put (9,65) {-3}
\put (39,65){-2}
\put (69,65){-1}
\put (10,50){0}
\put (40, 50){1}
\put (70,50){2}
\put (10,35){3}
\put (40,35){4}
\put (70,35){5}
\put (10,20){6}
\put (40,20){7}
\put (70,20){8}
\put (10,5){9}
\put (38,5){10}
\put (68,5){11}
\put (72.5,53){\circle{13}}
\put (72.5,68){\circle{13}}
\put (72.5,38){\circle{13}}
\put (12.5,68){\circle{13}}
\put (12.5,8){\circle{13}}
\put (72.5,83){\circle{13}}
\put (42.5,83){\circle{13}}
\put (12.5,83){\circle{13}}
\put (42.5,-3){.}
\put (42.5,-7){.}
\put (42.5,-11){.}
\put (72.5,-3){.}
\put (72.5,-7){.}
\put (72.5,-11){.}
\put (11.5,-3){.}
\put (11.5,-7){.}
\put (11.5,-11){.}
\end{picture}
\end{array}
$$
\end{example}

\begin{remark}
Note that an abacus for $\lambda$ is not unique because it depends on the set of $\beta$-numbers chosen for $\lambda$.  However, from any abacus of $\lambda$ one can obtain the partition $\lambda$ by counting the number of gaps before every bead in the abacus in reading order. In the first example above for instance, we see that $\lambda_1 = 8$ since the eight numbers 10,9,8,6,5,3,2,0 are exactly the eight gaps before the bead corresponding to the final bead at position $11$. We will say that a bead is \textit{active} if it 
occurs in the positions between the first gap and the last bead, in reading order.  The active beads are those that correspond to a nonzero part in a partition.  In the left example above, the bead in spot 11 is active since it corresponds to the part $\lambda_1 = 8$, whereas the bead in spot -1 is not active since it corresponds to $\lambda_5 = 0$.
\end{remark}

\begin{definition}
We define the \textit{balance number} of an abacus to be the sum over all runners of the largest level in that runner which contains a bead. We say that an abacus is \textit{balanced} if its balance number is zero.
\end{definition}

\begin{example}
In the example above, the balance number of the first diagram is $-1+2+3 = 4$. The balance number for the second diagram is $3 + -2 + 1 = 2$, so neither are balanced.
\end{example}
\begin{remark}\label{r:unique_abacus}
Note that there is a unique abacus which represents a given partition for each balance number. In particular, there is a unique abacus of $\lambda$ with balance number 0.  The balance number for a set of $\beta$-numbers of $\lambda$ will increase by exactly $1$ when the vector \textbf{1} is added to the set of $\beta$-numbers. On the abacus picture, this corresponds to shifting all of the beads forward one entry in the reading order.
\end{remark}

\begin{definition}
A runner is called \textit{flush} if no bead on the runner is preceded in reading order by a gap on that same runner.  We say that an abacus is \textit{flush} if every runner is flush.
\end{definition}

\begin{theorem}\cite[Theorem 2.7.16, Lemma 2.7.38]{JK}\label{t:flush_abacus}
$\lambda$ is an $\ell$-core if and only if any (equivalently, every) abacus of $\lambda$ is flush.  Moreover, in the balanced flush abacus of an $\ell$-core $\lambda$, each active bead on runner $i$ corresponds to a row of $\lambda$ whose rightmost box has residue $i$.
\end{theorem}

In the case that the corresponding abacus is not balanced, the boxes corresponding to the active beads on runner $i$ will share the same residue, but the residue may not be $i$.

\begin{example}\label{52111}
One can check that the partition $\lambda = (5,2,1,1,1)$ is a 4-core. One set of $\beta$-numbers for $\lambda$ is $(8,4,2,1,0,-2,-3,-4, \dots )$. This abacus is balanced as $2+0+0+(-2)=0$. All of the runners are flush. The active beads on runner 0 lie in positions 8, 4, 0 and these correspond to rows 1, 2 and 5 of the partition diagram whose final box of residue 0 is highlighted.

\begin{figure}[ht]
\begin{tabular}{llll}
\begin{picture}(105,40)(0,80)
\put (42.5,96){.}
\put (42.5,92){.}
\put (42.5,100){.}
\put (72.5,96){.}
\put (72.5,92){.}
\put (72.5,100){.}
\put (102.5,96){.}
\put (102.5,92){.}
\put (102.5,100){.}
\put (11.5,96){.}
\put (11.5,92){.}
\put (11.5,100){.}
\put (9,80) {-8}
\put (39,80){-7}
\put (69,80){-6}
\put (99,80){-5}
\put (9,65) {-4}
\put (39,65){-3}
\put (69,65){-2}
\put (99,65){-1}

\put (10,50){\bf 0}
\put (40, 50){1}
\put (70,50){2}
\put (101,50){3}

\put (10,35){\bf 4}
\put (40,35){5}
\put (70,35){6}
\put (101,35){7}

\put (10,20){\bf 8}
\put (40,20){9}
\put (68,20){10}
\put (99,20){11}

\put (72.5,53){\circle{13}}
\put (72.5,68){\circle{13}}
\put (72.5,83){\circle{13}}

\put (102.5,83){\circle{13}}

\put (42.5,53){\circle{13}}
\put (42.5,68){\circle{13}}
\put (42.5,83){\circle{13}}

\put (12.5,23){\circle{13}}
\put (12.5,38){\circle{13}}
\put (12.5,53){\circle{13}}
\put (12.5,68){\circle{13}}
\put (12.5,83){\circle{13}}

\put (42.5,12){.}
\put (42.5,8){.}
\put (42.5,4){.}
\put (72.5,12){.}
\put (72.5,8){.}
\put (72.5,4){.}
\put (102.5,12){.}
\put (102.5,8){.}
\put (102.5,4){.}
\put (11.5,12){.}
\put (11.5,8){.}
\put (11.5,4){.}
\end{picture}
& & \hspace{0.5in}  &
\parbox[t]{2in}{
\tableau{\mbox{0} & \mbox{1} & \mbox{2} & \mbox{3} & \mbox{\bf 0} \\
         \mbox{3} & \mbox{\bf 0} \\
         \mbox{2} \\
         \mbox{1} \\
         \mbox{\bf 0} \\
} } \\
 & & & \\
\end{tabular}
\caption{This abacus represents the 4-core $(5,2,1,1,1)$.  The boxes of the corresponding partition diagram have been filled with their residue. }\label{f:ex_pi}
\end{figure}
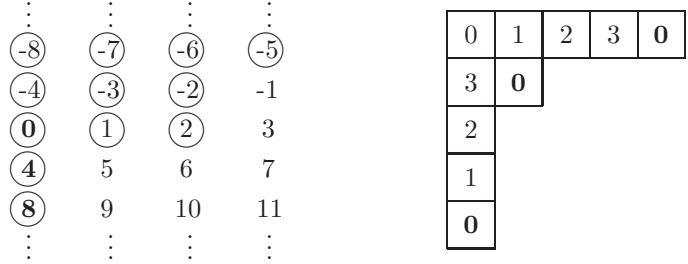
\end{example}

\subsection{The bijection on abacus configurations}

Here we describe the bijection $\Phi_{\ell}^k : \mathcal{C}_{\ell}^k \rightarrow \mathcal{C}_{\ell-1}^{\leq k}$.  Given $\lambda \in \mathcal{C}_{\ell}^k$ and an abacus for $\lambda$, remove the whole runner which contains the largest bead. Place the remaining runners into an $(\ell-1)$ abacus in order. In other words, renumber the runners $0, \ldots, \ell-2$, keeping the levels of the entries as before. This will correspond to an $(\ell - 1)$-core $\mu$ with largest part at most $k$. Then we define $\Phi_{\ell}^k : \mathcal{C}_{\ell}^k \rightarrow \mathcal{C}_{\ell-1}^{\leq k}$ to be the map which takes $\lambda$ to $\mu$. Observe that $ \Phi_{\ell}^k$ is well defined, independent of the choice of abacus for $\lambda$. 

To see that $\Phi_{\ell}^k$ is a bijection, observe that the map can be reversed.  Starting from an abacus of the $(\ell-1)$-core $\mu$, insert a new flush runner whose largest bead occurs just after the $k^{\textrm{th}}$ gap in the reading order.  This yields a flush abacus for the $\ell$-core $\lambda$ with $\lambda_1 = k$.

\begin{example}
Let $\ell = 4$ and $\lambda = (8,5,2,2,1,1,1)$. An abacus for $\lambda$ is:

\begin{center}
\begin{picture}(80,100)
\put (42.5,96){.}
\put (42.5,92){.}
\put (42.5,100){.}
\put (72.5,96){.}
\put (72.5,92){.}
\put (72.5,100){.}
\put (102.5,96){.}
\put (102.5,92){.}
\put (102.5,100){.}
\put (11.5,96){.}
\put (11.5,92){.}
\put (11.5,100){.}
\put (9,80) {-8}
\put (39,80){-7}
\put (69,80){-6}
\put (99,80){-5}
\put (9,65) {-4}
\put (39,65){-3}
\put (69,65){-2}
\put (99,65){-1}

\put (10,50){0}
\put (40, 50){1}
\put (70,50){2}
\put (101,50){3}

\put (10,35){4}
\put (40,35){5}
\put (70,35){6}
\put (101,35){7}

\put (10,20){8}
\put (40,20){9}
\put (68,20){10}
\put (99,20){11}

\put (72.5,53){\circle{13}}
\put (72.5,68){\circle{13}}
\put (72.5,83){\circle{13}}

\put (102.5,83){\circle{13}}

\put (42.5,53){\circle{13}}
\put (42.5,68){\circle{13}}
\put (42.5,83){\circle{13}}

\put (72.5,23){\circle{13}}
\put (72.5,38){\circle{13}}
\put (102.5,68){\circle{13}}
\put (12.5,83){\circle{13}}

\put (42.5,12){.}
\put (42.5,8){.}
\put (42.5,4){.}
\put (72.5,12){.}
\put (72.5,8){.}
\put (72.5,4){.}
\put (102.5,12){.}
\put (102.5,8){.}
\put (102.5,4){.}
\put (11.5,12){.}
\put (11.5,8){.}
\put (11.5,4){.}
\end{picture}
\end{center}

The largest $\beta$-number is 10. Removing the whole runner containing the 10, we get the remaining diagram with runners relabeled for $\ell = 3$
\begin{center}
\begin{picture}(80,100)
\put (42.5,96){.}
\put (42.5,92){.}
\put (42.5,100){.}
\put (72.5,96){.}
\put (72.5,92){.}
\put (72.5,100){.}
\put (102.5,96){.}
\put (102.5,92){.}
\put (102.5,100){.}
\put (11.5,96){.}
\put (11.5,92){.}
\put (11.5,100){.}
\put (9,80) {-6}
\put (39,80){-5}
\put (69,80){$\times$}
\put (99,80){-4}
\put (9,65) {-3}
\put (39,65){-2}
\put (69,65){$\times$}
\put (99,65){-1}
\put (10,50){0}
\put (40, 50){1}
\put (69,50){$\times$}
\put (101,50){2}
\put (10,35){3}
\put (40,35){4}
\put (69,35){$\times$}
\put (101,35){5}
\put (10,20){6}
\put (40,20){7}
\put (69,20){$\times$}
\put (101,20){8}
\put (72.5,53){\circle{13}}
\put (72.5,68){\circle{13}}
\put (72.5,83){\circle{13}}
\put (102.5,83){\circle{13}}
\put (42.5,83){\circle{13}}
\put (42.5,68){\circle{13}}
\put (42.5,53){\circle{13}}
\put (72.5,22){\circle{13}}
\put (72.5,38){\circle{13}}
\put (102.5,68){\circle{13}}
\put (12.5,83){\circle{13}}
\put (42.5,12){.}
\put (42.5,8){.}
\put (42.5,4){.}
\put (72.5,12){.}
\put (72.5,8){.}
\put (72.5,4){.}
\put (102.5,12){.}
\put (102.5,8){.}
\put (102.5,4){.}
\put (11.5,12){.}
\put (11.5,8){.}
\put (11.5,4){.}
\end{picture}
\end{center}
These are a set of $\beta$-numbers for the partition $(2,1,1)$, which is a 3-core with largest part $\leq 8$.  For the reverse bijection when $k=8$, notice that the eighth gap is at entry 7 which dictates where we insert the new runner and beads.  Also note in this example that the first abacus has balance number -1 while its image has balance number -3, so balance number is not necessarily preserved.
\end{example}

\subsection{The bijection on core partitions}\label{s:phi_on_diagram}

Another way to describe $\Phi_{\ell}^k$ is on the Young diagram of $\lambda$. Applying $\Phi_{\ell}^k$ to $\lambda$ is the same as removing all of the rows $i$ of $\lambda$ for which $ h_{(i,1)} \equiv h_{(1,1)} \mod \ell$. 
To illustrate, we show the bijection on the same example $\lambda = (8,5,2,2,1,1,1)$, but performed on a Young diagram instead of an abacus. We start by drawing the Young diagram and writing the hooks lengths of the boxes in the first column. The bijection simply deletes the rows which have a hook length in the first column equivalent to the hook length $h_{(1,1)}^{\lambda}\mod \ell$.
\noindent

$$
\tableau{14 & \mbox{} & \mbox{} & \mbox{} & \mbox{} & \mbox{} & \mbox{} & \mbox{} \\
10 & & \mbox{} & \mbox{} & \mbox{} \\
6 & \mbox{}\\
{\bf 5} & \mbox{}\\
{\bf 3}\\
2\\
{\bf 1}} 
\put (-254,8){\small Hook length $\equiv 14 \, mod \, 4 \to$}
\put (-254,-10){\small Hook length $\equiv 14 \, mod \, 4 \to$}
\put (-254,-28){\small Hook length $\equiv 14 \, mod \, 4 \to$}
\put (-254,-82){\small Hook length $\equiv 14 \, mod \, 4 \to$}
\put (10,-65){ \vspace{0.2in} $\stackrel{\Phi_{4}^8}{\mapsto}$}
\put (40,-71){\tableau{\mbox{} & \mbox{} \\ \mbox{} \\ \mbox{}}}
$$

Deleting the corresponding rows, we get that $\Phi_4^8 (8,5,2,2,1,1,1) = (2,1,1)$.
\normalsize

\section{Cores and the action of $\widetilde{S_{\ell}}$ on the finite root lattice}\label{s:geometry}

In this section we recall that the $\ell$-cores index a system of minimal length coset representatives for $\widetilde{S_{\ell}} / S_{\ell}$ and describe some associated geometry.

\subsection{The affine root system}\label{s:aff_geometry}
Following \cite{humphreys}, let $\{ \e_1, \e_2, \dots, \e_{\ell} \}$ be an orthonormal basis of the Euclidean space $\mathbf{R}^{\ell}$ and denote the corresponding inner product by $(\cdot, \cdot)$.
For $1 \leq i \leq \ell - 1$, let $s_i$ be the reflection defined by interchanging $\e_i$ and $\e_{i+1}$; the reflecting hyperplanes are discussed below.  Then $\{s_1, \ldots, s_{\ell-1}\}$ are a set of Coxeter generators for the symmetric group $S_{\ell}$, which acts on $\mathbf{R}^{\ell}$ by permuting coordinates in the $\e_i$ basis.

Let $s_0$ be the affine reflection of $\mathbf{R}^{\ell}$ defined on $v = \sum_{j=1}^{\ell} a_j \e_j$ by
\[ s_0(v) = (a_{\ell} + 1) \e_1 + a_2 \e_2 + \cdots + a_{\ell-1} \e_{\ell-1} + (a_{1} - 1) \e_{\ell}. \]
Define the \em simple roots \em $\Delta$ \em of type $A_{\ell-1}$ \em to be the collection of $\ell-1$ vectors
\[ \a_1 = \e_1 - \e_2, \ \ \a_2 = \e_2 - \e_3, \ \ \ldots, \ \ \a_{\ell-1} = \e_{\ell-1} - \e_{\ell}. \]

The $\mathbf{Z}$-span $\Lambda_R$ of $\Delta$ is called the \em root lattice of type $A_{\ell-1}$\em.
Let $V = \mathbf{R} \otimes_{\mathbf{Z}} \Lambda_R \subsetneq \mathbf{R}^{\ell}$.  Observe that each reflection $s_i$ preserves $V$ and so
$\{s_0, s_1, \ldots, s_{\ell-1}\}$ are a set of Coxeter generators for the affine symmetric group $\widetilde{S_{\ell}}$ acting on $V$.  For $w \in \widetilde{S_{\ell}}$ we let $l(w)$ denote Coxeter length.  From now on, we restrict our attention from $\mathbf{R}^{\ell}$ to $V$.  

In this presentation we see that $S_{\ell}$ is a parabolic subgroup of $\widetilde{S_{\ell}}$.  We form the parabolic quotient 
\[ \widetilde{S_{\ell}} / S_{\ell} = \{ w \in \widetilde{S_{\ell}} : l(w s_i) > l(w) \text{ for all $s_i$ where $1 \leq i \leq \ell-1$} \}. \]
By a standard result in the theory of Coxeter groups, this set gives a unique representative of minimal length from each coset $w S_{\ell}$ of $\widetilde{S_{\ell}} / S_{\ell}$.  For more on this construction, see \cite[Section 2.4]{bjorner-brenti}.
Another standard result is that $\widetilde{S_{\ell}}$ acts on $V$ as the semidirect product of $S_{\ell}$ and the translation group corresponding to the root lattice $\Lambda_R$.  
Hence, $\widetilde{S_{\ell}} / S_{\ell}$ is also in bijection with $\Lambda_R$ and we identify $\Lambda_R$ with the translation subgroup $\{t_{\mathbf{a}} : \mathbf{a} \in \Lambda_R\}$ of $\widetilde{S_{\ell}}$.

Let us consider this situation geometrically.  Denote the set of \em finite roots \em by $\Pi = \{ w \a_{i} : w \in S_{\ell}, \a_i \in \Delta \} \subset V$.  It is a standard fact that each root $\a \in \Pi$ can be written as an integral linear combination of the simple roots $\Delta$ such that all of the coefficients are positive or all coefficients are negative.  Therefore, $\Pi$ can be decomposed as $\Pi = \Pi^{+} \uplus \Pi^{-}$.

For each finite root $\a$ and integer $k$ we can define an affine hyperplane
\[ H_{\a, k} = \{ v \in V : (v, \a) = k \}. \]
Observe that $s_i$ is the reflection over the hyperplane $H_{\a_i,0}$ for $1 \leq i < \ell$ while $s_0$ is the reflection over $H_{\theta, 1}$ where $\theta = \e_1 - \e_{\ell} = \sum_{i=1}^{\ell-1} \a_i$.  Let $\mathcal{H}$ denote the collection of all affine hyperplanes $H_{\a, k}$ for $\a \in \Pi, k \in \mathbf{Z}$.  Let $\mathcal{A}$ be the set of all connected components of $V \setminus \bigcup_{H \in \mathcal{H}} H$.  Each element of $\mathcal{A}$ is called an \em alcove\em.  In particular,
\[ A_{\circ} = \{ v \in V : 0 < (v, \a) < 1 \text{ for all } \a \in \Pi^{+} \} \]
is called the \em fundamental alcove \em whose closure is a fundamental domain for the action of $\widetilde{S_{\ell}}$ on $V$.  

\begin{proposition}\cite[Section 4.5]{humphreys}
The affine Weyl group $\widetilde{S_{\ell}}$ permutes the collection of alcoves transitively and freely.
The closure of $A_{\circ}$ is a fundamental domain for the action of $\widetilde{S_{\ell}}$ on $V$.
\end{proposition}

Define 
\[ B_{\circ} = \bigcup_{w \in S_{\ell}} w A_{\circ}. \]
The set $B_{\circ}$ contains one alcove for each permutation in $S_{\ell}$ and the closure of $B_{\circ}$ is a fundamental domain for the action of translation by $\Lambda_R$ on $V$.  Moreover, for any affine permutation $t_{\mathbf{a}} w \in \Lambda_R \rtimes S_{\ell}$ we have that $B_{\circ} + \mathbf{a}$ corresponds to the left coset of $S_{\ell}$ containing $t_{\mathbf{a}} w$ in $\widetilde{S_{\ell}} / S_{\ell}$.
In fact, the set of all $B_{\circ} + \mathbf{a}$ is precisely the $\widetilde{S_{\ell}}$-orbit of $B_{\circ}$.
Since length can be computed by counting the minimal number of hyperplanes which must be crossed in a path back to $A_{\circ}$, we have that $\mathbf{a}$ determines the alcove of $B_{\circ} + \mathbf{a}$ that represents $t_{\mathbf{a}} w$ and has minimal length:  it is the alcove which requires the fewest such hyperplane crossings.  This describes a bijection that we denote $\varpi : \Lambda_R \rightarrow \widetilde{S_{\ell}} / S_{\ell}$ in which $\mathbf{a} \mapsto t_{\mathbf{a}} w_{\mathbf{a}}$.
Figure~\ref{f:sl3} shows the Euclidean space $V$ associated to type $A_{2}$ in the context of the bijection $\Phi_{\ell}^k$.

\subsection{$\ell$-cores are minimal length coset representatives}

We now show how an $\ell$-core can be associated to each $\mathbf{a} \in \Lambda_R$.
Let $\mathbf{a} = (a_1, \ldots, a_{\ell})$ be a vector in $\Lambda_R$ written with respect to the $\e_i$ basis, so each $a_i \in \mathbf{Z}$ and $\sum_{i = 1}^{\ell} a_i = 0$.  We form a balanced flush abacus from $\mathbf{a}$ by filling the $(i-1)^{\mathrm{st}}$ runner with beads from $-\infty$ down to level $a_i$.  By Remark~\ref{r:unique_abacus} and Theorem~\ref{t:flush_abacus}, every $\ell$-core has exactly one balanced flush abacus.  Hence, we obtain bijections whose composition we denote by $\pi$.
\[ \pi: \{ (a_1, \ldots, a_{\ell}) : a_i \in \mathbf{Z}, \sum_{i=1}^{\ell} a_i = 0 \} \rightarrow \{ \text{balanced flush abaci} \} \rightarrow \mathcal{C}_{\ell}. \]

\begin{remark}
Because the runners of an abacus are usually labeled by $0 \leq i \leq \ell-1$ but coordinates of $\mathbf{R}^{\ell}$ are labeled by $1 \leq j \leq \ell$ we will sometimes coordinatize $\Lambda_R$ as $\mathbf{a} = (a_1, \ldots, a_{\ell})$ or as $\mathbf{b} = (b_0, \ldots, b_{\ell-1})$.
\end{remark}

\begin{example}
Let $\ell = 4$ and let $\mathbf{a} = 2\e_1 + 0\e_2 +0\e_3 -2 \e_4$. Then we draw an abacus as shown in Figure~\ref{f:ex_pi} above with beads down to level 2 in runner 0, level 0 in runner 1, level 0 in runner 2, and level -2 in runner 3.
\end{example}

Next, we observe that the $\ell$-cores inherit an action of $\widetilde{S_{\ell}}$ from the bijection $\pi$.  To describe the action, we draw the diagram of $\lambda = (\lambda_1, \ldots, \lambda_r)$ and fill each box $(i,j)$ with the residue $j-i \mod \ell$.  Fix $\lambda \in \mathcal{C}_{\ell}$ and suppose $\mathbf{b} = (b_0, \ldots, b_{\ell-1}) = \pi^{-1}(\lambda)$.  For $0 \leq i \leq \ell-1$, we say that $s_i$ is an \em ascent \em for $\lambda$ if $b_{i-1} > b_{i} - \delta_{i,0}$, and we say that $s_i$ is a \em descent \em for $\lambda$ if $b_{i-1} < b_{i} - \delta_{i,0}$.  Here, we interpret $b_{-1}$ as $b_{\ell-1}$.

\begin{remark}
Note that the definition for $s_i$ to be a descent (respectively, ascent, neither) given above corresponds to $l(\varpi(s_i \mathbf{b})) < l(\varpi(\mathbf{b}))$ (respectively, $>$, $=$).  This corresponds to whether $l(s_i \varpi(\mathbf{b})) < l(\varpi(\mathbf{b}))$.
\end{remark}

\begin{example}
The 4-core $(5,2,1,1,1)$ corresponds to $w = s_0 s_1 s_2 s_3 s_2 s_1 s_0 = t_{(2,0,0,-2)}$.  Hence, $s_1$ and $s_3$ are ascents for $\lambda$ and $s_0$ is a descent for $\lambda$.  Observe that $s_2$ is neither an ascent nor a descent because $s_2 w = s_0 s_1 s_2 s_3 s_2 s_1 s_0 s_2$ ceases to be a minimal length left coset representative.  Correspondingly $s_2 (2, 0, 0, -2) = (2, 0, 0, -2)$.
\end{example}

\begin{proposition}\label{p:action}
Let $\lambda$ be an $\ell$-core.  If $s_i$ is an ascent for $\lambda$ then $s_i$ acts on $\lambda$ by adding all boxes with residue $i$ to $\lambda$ such that the result is a partition.  If $s_i$ is a descent for $\lambda$ then $s_i$ acts on $\lambda$ by removing all of the boxes with residue $i$ that lie at the end of both their row and column so that their removal results in a partition.  If $s_i$ is neither an ascent nor a descent for $\lambda$ then $s_i$ does not change $\lambda$.
\end{proposition}

\begin{proof}
Begin by considering the action of $s_i$ on abaci that comes from the action on $\Lambda_R \subset V$.  Applying $s_i$ for $1 \leq i \leq \ell-1$ corresponds to exchanging adjacent runners $i-1$ and $i$ in the balanced flush abacus whose runners are labeled $0, \ldots, \ell-1$.  Applying the $s_0$ generator first adds a bead to runner $\ell-1$ and removes a bead from runner $0$ so that they stay flush, and then exchanges the two runners.  

Since the coordinates of $\mathbf{b} \in \Lambda_R$ sum to 0, Theorem~\ref{t:flush_abacus} implies that each active bead in runner $i$ of the balanced flush abacus corresponds to a row of $\lambda$ whose rightmost box has residue $i$.  Because the abacus is flush, exchanging runners $i$ and $i-1$ either adds some set of boxes with residue $i$ to the diagram of $\lambda$ in the case that $s_i$ is an ascent, or else removes a set of boxes with residue $i$ in the case that $s_i$ is a descent.  The result is again a balanced flush abacus so corresponds to an $\ell$-core.  If $s_i$ is neither an ascent nor a descent then $b_i = b_{i-1} + \delta_{i,0}$ so the abacus remains unchanged.

Observe that a box with residue $i$ is removable if and only if it lies at the end of its row and column.  This occurs if and only if it corresponds to an active bead on runner $i$ with a gap immediately preceding it in the reading order of the abacus.  Similarly, a box with residue $i$ is addable if and only if it corresponds to a gap on runner $i$ with an active bead immediately succeeding it in the reading order of the abacus.  The action of $s_i$ swaps runners $i$ and $i-1$ which therefore interchanges all of the $i$-addable and $i$-removable boxes.
\end{proof}

\begin{remark}
This action can be described as a special case of the action on the crystal graph associated to the irreducible highest-weight representation $V(\Lambda_0)$ of $\widehat{\mathfrak{sl}_{\ell}}$ that is given by the operators $\widetilde{f_i}^{\varphi_i - \varepsilon_i}$ or $\widetilde{e_i}^{\varepsilon_i - \varphi_i}$.  For details, see \cite{G}, \cite{Kl} or \cite{MM}.
\end{remark}

Let $\lambda$ be an $\ell$-core.  Then we can recursively define a canonical reduced expression for $\varpi(\pi^{-1}(\lambda))$ that we denote $w(\lambda)$ by choosing $w(\lambda) = s_i w(\widehat{\lambda})$ where $i$ is the residue of the rightmost box in the bottom row of $\lambda$ and $\widehat{\lambda}$ is the result of applying $s_i$ to $\lambda$ as in Proposition~\ref{p:action}.  Note that $s_i$ is always a descent for this choice of $i$.  The empty partition corresponds to the identity Coxeter element.  This reduced expression was previously defined in \cite[Definition 45]{LM}.  

\begin{example}
The canonical reduced expression for the 4-core $(5, 2, 1, 1, 1)$ shown in Figure~\ref{f:ex_pi} is $s_0 s_1 s_2 s_3 s_2 s_1 s_0$.  The first step in reducing this expression to the identity removes the three boxes labeled 0 that lie at the end of their rows which is recorded as the leftmost $s_0$ in the expression.
\[
\tableau{\mbox{0} & \mbox{1} & \mbox{2} & \mbox{3} & \mbox{\bf 0} \\
         \mbox{3} & \mbox{\bf 0} \\
         \mbox{2} \\
         \mbox{1} \\
         \mbox{\bf 0} \\
} \stackrel{s_0}{\hspace{0.2in} \longleftarrow \hspace{0.2in}}
\tableau{\mbox{0} & \mbox{1} & \mbox{2} & \mbox{3} \\
         \mbox{3} \\
         \mbox{2} \\
         \mbox{\bf 1} \\
} \stackrel{s_1}{\hspace{0.2in} \longleftarrow \hspace{0.2in}}
\tableau{\mbox{0} & \mbox{1} & \mbox{2} & \mbox{3} \\
         \mbox{3} \\
         \mbox{\bf 2} \\
} 
\]
\[
\stackrel{s_2}{\hspace{0.2in} \longleftarrow \hspace{0.2in}}
\tableau{\mbox{0} & \mbox{1} & \mbox{2} & \mbox{\bf 3} \\
         \mbox{\bf 3} \\
} \stackrel{s_3}{\hspace{0.2in} \longleftarrow \hspace{0.2in}}
\tableau{\mbox{0} & \mbox{1} & \mbox{\bf 2} \\
} \stackrel{s_2}{\hspace{0.2in} \longleftarrow \hspace{0.2in}}
\tableau{\mbox{0} & \mbox{\bf 1} \\
} 
\]
\[
\stackrel{s_1}{\hspace{0.2in} \longleftarrow \hspace{0.2in}}
\tableau{\mbox{\bf 0} \\
} \stackrel{s_0}{\hspace{0.2in} \longleftarrow \hspace{0.2in}}
\emptyset.
\]
\end{example}

\begin{proposition}\label{p:redexp}
For $\lambda \in \mathcal{C}_{\ell}$ we have that $w(\lambda)$ is a reduced expression for the minimal length coset representative indexed by $\lambda$.

The Coxeter length of $w(\lambda)$ is 
\[ l(w(\lambda)) = \sum_{i = 0}^{\ell-1} \lambda_{R(i)} \]
where $R(i)$ is the longest row of $\lambda$ whose rightmost box has residue $i$.
\end{proposition}

\begin{proof}
Let $i$ be the residue of the rightmost box in the bottom row $\lambda_r$ of $\lambda$.  Then the balanced flush abacus configuration corresponding to $\lambda$ has an active bead $B$ representing $\lambda_r$ on runner $i$ by Theorem~\ref{t:flush_abacus}.  Observe that $B$ has a gap immediately preceding it in the reading order because $B$ is the first active bead in the reading order and $\lambda_r \neq 0$.  Since the abacus is flush, every box of $\lambda$ with residue $i$ that lies at the end of its row corresponds to some active bead $B'$ on runner $i$ of the abacus with a gap immediately preceding $B'$ in the reading order.  Hence, every box of $\lambda$ with residue $i$ that lies at the end of its row also lies at the end of its column and applying $s_i$ removes every such box by Proposition~\ref{p:action}.

Iterating this process eventually produces the empty partition, corresponding to the identity Coxeter element.  At each step, we have shown that applying $s_i$ removes exactly one box from $R(i)$ yielding the length formula.
\end{proof}

\begin{example} Let $\lambda = (10,7,4,3,2,2,2,1,1,1)$ and $\ell = 4$. Then the biggest part ending in residue 1 is $\lambda_1 = 10$, the biggest part ending in residue 3 is $\lambda_4 = 3$, the biggest part ending in residue 0 is $\lambda_6 = 2$, and no part ends in residue 2. Hence the Coxeter length of the minimal length coset representative for this 4-core is $10+3+2 = 15$. The canonical minimal length coset representative for this core is $w(\lambda) = s_3s_0s_1s_2s_3s_0s_1s_3s_2s_1s_0s_3s_2s_1s_0$.
$$\tableau{0&1&2&3&0&1&2&3&0&1\\
3&0&1&2&3&0&1\\
2&3&0&1\\
1&2&3\\
0&1\\
3&0\\
2&3\\
1\\
0\\
3}$$
\end{example}

There is another way to obtain the root vector $\pi^{-1}(\lambda) = (b_0, \ldots, b_{\ell-1}) \in \Lambda_R$ from the partition $\lambda \in \mathcal{C}_{\ell}$ due to Garvan, Kim and Stanton \cite{garvan-kim-stanton}.  Say that \em region $r$ \em of the diagram of $\lambda$ is the set of boxes $(i,j)$ satisfying $(r-1) \ell \leq j-i < r \ell$.  We call a box \em row-exposed \em if it lies at the end of its row.  Then, set $b_i$ to be the maximum region of $\lambda$ which contains a row-exposed box with residue $i$ for $0 \leq i \leq \ell-1$.  In particular, we pad $\lambda$ with parts of size zero if necessary and label all of the boxes before the $0^{\textrm{th}}$ column by their residue.  In this way $b_i$ is well-defined because column 0 contains infinitely many row-exposed boxes.  We call the vector $(b_0, \ldots, b_{\ell-1})$ obtained in this fashion the \em $n$-vector of $\lambda$ \em and we show that it is the same vector as $\pi^{-1}(\lambda)$.

\begin{example}
Let $\ell = 4$ and $\lambda = (6,3,1,1)$. From the picture below, we see that the $n$-vector for $\lambda$ is $(-1,2,0,-1)$. 
$$
\tableau{0&1&2&3&0&1\\
3&0&1\\
2\\
1}
\put (-118,6){3}
\put (-118,-12){2}
\put (-118,-30){1}
\put (-118,-48){0}
\put (-118,-66){3}
\put (-118,-84){2}
\put (-118,-102){1}
\put (-118,-120){0}
\put (-30,-50){Region 1}
\put (-75,-85){Region 0}
\put (-100,-130){Region -1}
\put (10,-20){Region 2}
\put (-108,-124){\line(0,1){76}}
\put (10,6){$\leftarrow$ First 1}
\put (-84,-30){$\leftarrow$ First 2}
\put (-162,-66){First 3 $\to$}
\put (-162,-120){First 0 $\to$}
\linethickness{1.5pt}
\put (-36,18){\line(-1,0){18}}
\put (-36,18){\line(0,-1){18}}
\put (-36,0){\line(1,0){18}}
\put (-18,0){\line(0,-1){18}}
\put (-18,-18){\line(1,0){18}}
\put (0,-18){\line(0,-1){18}}
\put (0,-36){\line(1,0){18}}
\put (18,-36){\line(0,-1){18}}
\put (18,-54){\line(1,0){18}}
\put (36,-54){\line(0,-1){18}}
\put (36,-72){\line(1,0){18}}
\put (54,-72){\line(0,-1){18}}
\put (54,-90){\line(1,0){18}}
\put (72,-90){\line(0,-1){18}}
\put (72,-108){\line(1,0){18}}
\put (90,-108){\line(0,-1){18}}
\put (90,-126){\line(1,0){18}}
\put (108,-126){\line(0,-1){18}}
\put (108,-144){\line(1,0){18}}
\put (-108,18){\line(-1,0){18}}
\put (-108,18){\line(0,-1){18}}
\put (-108,0){\line(1,0){18}}
\put (-90,0){\line(0,-1){18}}
\put (-90,-18){\line(1,0){18}}
\put (-72,-18){\line(0,-1){18}}
\put (-72,-36){\line(1,0){18}}
\put (-54,-36){\line(0,-1){18}}
\put (-54,-54){\line(1,0){18}}
\put (-36,-54){\line(0,-1){18}}
\put (-36,-72){\line(1,0){18}}
\put (-18,-72){\line(0,-1){18}}
\put (-18,-90){\line(1,0){18}}
\put (0,-90){\line(0,-1){18}}
\put (0,-108){\line(1,0){18}}
\put (18,-108){\line(0,-1){18}}
\put (18,-126){\line(1,0){18}}
\put (36,-126){\line(0,-1){18}}
\put (36,-144){\line(1,0){18}}
\put (-108,-54){\line(-1,0){18}}
\put (-108,-54){\line(0,-1){18}}
\put (-108,-72){\line(1,0){18}}
\put (-90,-72){\line(0,-1){18}}
\put (-90,-90){\line(1,0){18}}
\put (-72,-90){\line(0,-1){18}}
\put (-72,-108){\line(1,0){18}}
\put (-54,-108){\line(0,-1){18}}
\put (-54,-126){\line(1,0){18}}
\put (-36,-126){\line(0,-1){18}}
\put (-36,-144){\line(1,0){18}}
\put (-108,-126){\line(-1,0){18}}
\put (-108,-126){\line(0,-1){18}}
\put (-108,-144){\line(1,0){18}}
$$
\end{example}

\begin{lemma}\label{l:n_vector}
Let $(b_0, \ldots, b_{\ell-1})$ be the $n$-vector of an $\ell$-core $\lambda$ and let $s_i(\lambda)$ be the result of applying $s_i$ to $\lambda$ as described in Proposition~\ref{p:action}.  Then, the $n$-vector of $s_i(\lambda)$ is $(b_0, \dots, b_{i}, b_{i-1}, \dots b_{\ell-1})$ for $1 \leq i \leq \ell-1$, or $(b_{\ell-1} + 1, b_2, \dots, b_{\ell-2}, b_{0}-1)$ if $i = 0$.
\end{lemma}
\begin{proof}
Denote the $n$-vector of $s_i(\lambda)$ by $(n_0, \dots, n_{\ell-1})$.  It follows from the proof in \cite[Bijection 2]{garvan-kim-stanton} that if there exists a row-exposed $i$-box in region $r$ then there exist row-exposed $i$-boxes in each region $< r$.  

Observe that a row $\lambda_p$ of $\lambda$ has an addable $i$-box in region $r$ exactly if:
\begin{enumerate}
\item the row $\lambda_p$ has a row-exposed $(i-1)$-box in region $r - \delta_{i,0} \leq b_{i-1}$, and
\item the row $\lambda_{p-1}$ does not have a row-exposed $i$-box, so $b_i < r$.
\end{enumerate}

Similarly, a row $\lambda_q$ has a removable $i$-box in region $s$ exactly if:
\begin{enumerate}
\item the row $\lambda_q$ has a row-exposed $i$-box in region $s \leq b_i$, and
\item the row $\lambda_{q+1}$ does not have a row-exposed $(i-1)$-box, so $b_{i-1} < s - \delta_{i,0}$.
\end{enumerate}

Moreover, an $i$-box is addable (removable) in $\lambda$ only if it is removable (addable, respectively) in $s_i(\lambda)$.  Therefore, when we apply $s_i$ to $\lambda$ we interchange each of the regions $r$ and $s$ satisfying
\[ b_i < r \leq b_{i-1} + \delta_{i,0} \]
and
\[ b_{i-1} + \delta_{i,0} < s \leq b_{i}. \]

Hence, we obtain
\[ n_i = b_{i-1} + \delta_{i,0} \]
and
\[ n_{i-1} = b_i - \delta_{i,0}. \]
\end{proof}

\begin{corollary}
The $n$-vector of $\lambda$ is $\pi^{-1}(\lambda)$.
\end{corollary}
\begin{proof}
The $n$-vector of the empty partition is equal to $\pi^{-1}(\emptyset) = (0, \ldots, 0)$.  The result then follows by induction on the Coxeter length of $w(\lambda)$ by Lemma~\ref{l:n_vector}.
\end{proof}

\begin{proposition}\label{p:k_from_vector}
Suppose that $\pi(\mathbf{a}) = \pi(a_1, \ldots, a_{\ell}) = \lambda$.  Then we have 
\[ \lambda_1 = (a_i - 1) \ell + i \]
where $a_i$ is the rightmost occurrence of the largest coordinate in $\mathbf{a}$.  Also, 
\[ \lambda_1 = \sum_{j = 1}^{i-1} (a_i - a_j) + \sum_{j=i+1}^{\ell} (a_i - a_j - 1). \]
\end{proposition}
\begin{proof}
Consider the balanced flush abacus corresponding to $\lambda$.  Then $\lambda_1$ corresponds to the last active bead $B$ in the reading order, and $B$ lies on runner $i-1$.  In particular, if there are multiple occurrences of the largest coordinate in $\mathbf{a}$ then $\lambda_1$ corresponds to the bead on the rightmost runner.  The number of boxes in $\lambda_1$ is the number of gaps prior to $B$ in the reading order of the abacus and the second formula follows from counting these gaps, using that the beads are flush on each runner.

Since the abacus is balanced, we have that the number of beads strictly below the zero level must be equal to the number of gaps weakly above the zero level.  If the last active bead occurs in level $j$ of runner $i-1$ then we could move all of the beads below the zero level to fill in the gaps above the zero level, and so count the gaps starting from entry 0 to the entry that contained $B$ as $(j-1) \ell + i$.  This yields the first formula.
\end{proof}

\begin{example}
The balanced flush abacus corresponding to the 4-core \[ \pi(1, -2, 2, -1) = \lambda = (7, 4, 3, 2, 1, 1, 1) \] is shown below together with the diagram in which the beads have been moved to calculate $\lambda_1 = 7$ as $(2-1) 4 + 3$.  Here, $a_3 = 2$ is the largest entry, corresponding to runner 2.

\begin{tabular}{ccc}
\begin{picture}(105,110)
\put (42.5,96){.}
\put (42.5,92){.}
\put (42.5,100){.}
\put (72.5,96){.}
\put (72.5,92){.}
\put (72.5,100){.}
\put (102.5,96){.}
\put (102.5,92){.}
\put (102.5,100){.}
\put (11.5,96){.}
\put (11.5,92){.}
\put (11.5,100){.}
\put (9,80) {-8}
\put (39,80){-7}
\put (69,80){-6}
\put (99,80){-5}
\put (9,65) {-4}
\put (39,65){-3}
\put (69,65){-2}
\put (99,65){-1}

\put (10,50){0}
\put (40, 50){1}
\put (70,50){2}
\put (101,50){3}

\put (10,35){4}
\put (40,35){5}
\put (70,35){6}
\put (101,35){7}

\put (10,20){8}
\put (40,20){9}
\put (68,20){10}
\put (99,20){11}

\put (5,45.5){\line(90,0){105}}

\put (72.5,53){\circle{13}}
\put (72.5,68){\circle{13}}
\put (72.5,83){\circle{13}}

\put (102.5,83){\circle{13}}

\put (42.5,83){\circle{13}}

\put (12.5,38){\circle{13}}
\put (12.5,53){\circle{13}}
\put (12.5,68){\circle{13}}
\put (12.5,83){\circle{13}}

\put (72.5,38){\circle{13}}
\put (72.5,23){\circle{13}}
\put (101.5,68){\circle{13}}

\put (42.5,12){.}
\put (42.5,8){.}
\put (42.5,4){.}
\put (72.5,12){.}
\put (72.5,8){.}
\put (72.5,4){.}
\put (102.5,12){.}
\put (102.5,8){.}
\put (102.5,4){.}
\put (11.5,12){.}
\put (11.5,8){.}
\put (11.5,4){.}
\end{picture} & \hspace{0.4in} &
\begin{picture}(105,100)
\put (42.5,96){.}
\put (42.5,92){.}
\put (42.5,100){.}
\put (72.5,96){.}
\put (72.5,92){.}
\put (72.5,100){.}
\put (102.5,96){.}
\put (102.5,92){.}
\put (102.5,100){.}
\put (11.5,96){.}
\put (11.5,92){.}
\put (11.5,100){.}
\put (9,80) {-8}
\put (39,80){-7}
\put (69,80){-6}
\put (99,80){-5}
\put (9,65) {-4}
\put (39,65){-3}
\put (69,65){-2}
\put (99,65){-1}

\put (10,50){0}
\put (40, 50){1}
\put (70,50){2}
\put (101,50){3}

\put (10,35){4}
\put (40,35){5}
\put (70,35){6}
\put (101,35){7}

\put (10,20){8}
\put (40,20){9}
\put (68,20){\underline{\bf{10}}}
\put (99,20){11}

\put (72.5,53){\circle{13}}
\put (72.5,68){\circle{13}}
\put (72.5,83){\circle{13}}

\put (102.5,83){\circle{13}}

\put (42.5,53){\circle{13}}
\put (42.5,68){\circle{13}}
\put (42.5,83){\circle{13}}

\put (12.5,53){\circle{13}}
\put (12.5,68){\circle{13}}
\put (12.5,83){\circle{13}}

\put (101.5,68){\circle{13}}
\put (103.5,53){\circle{13}}

\put (42.5,12){.}
\put (42.5,8){.}
\put (42.5,4){.}
\put (72.5,12){.}
\put (72.5,8){.}
\put (72.5,4){.}
\put (102.5,12){.}
\put (102.5,8){.}
\put (102.5,4){.}
\put (11.5,12){.}
\put (11.5,8){.}
\put (11.5,4){.}
\end{picture} \\
\end{tabular}
\end{example}

\begin{corollary}\label{c:hyperplane}
For $k \geq 0$, let $H_{\ell}^{k}$ denote the affine hyperplane
\[ H_{\ell}^{k} = \{ \mathbf{a} = (a_1, \ldots, a_{\ell}) \in \mathbf{R}^{\ell} : (\mathbf{a}, \e_{(k\ \mathrm{mod}\ \ell)}) = \lceil {k \over {\ell} } \rceil \} \cap V \]
inside $V$, where $1 \leq (k \mod \ell) \leq \ell$.
Then under the correspondence $\pi$, the $\ell$-cores $\lambda$ with $\lambda_1 = k$ all lie inside $H_{\ell}^{k} \bigcap \Lambda_R$.
\end{corollary}
\begin{proof}
We can write $k > 0$ uniquely as $(j - 1) \ell + i$ for $1 \leq i \leq \ell$ and $j \geq 1$.  In this case, $j = \lceil {k \over \ell} \rceil$ and $i \equiv k \mod \ell$.  The result then follows from the first formula of Proposition~\ref{p:k_from_vector}.  If $\lambda_1 = 0$ then $k = 0$ and $\lambda = \emptyset$ so the statement holds.
\end{proof}

\section{The bijection $\Phi_{\ell}^k$ is an affine linear isometry in $V$}\label{s:geometric_bijection}

\subsection{The bijection interpreted on $V$}

From Proposition~\ref{p:k_from_vector} we see that if $a_i$ is the rightmost occurrence of the largest coordinate in $\mathbf{a} = (a_1, \ldots, a_{\ell}) = \pi^{-1}(\lambda)$, then $\lambda_1 \equiv i \mod \ell$.  The next result describes $\Phi_{\ell}^{k}$ in terms of the root lattice coordinates.

\begin{theorem}\label{t:main}
Let $\psi_{\ell}$ be the affine map defined by $\psi_{\ell} (a_1, \ldots, a_{\ell}) = (a_{\ell}+1, a_1, a_2, \ldots, a_{\ell-1})$.  Then, 
\[ \pi^{-1} \circ \Phi_{\ell}^{k} \circ \pi(a_1, \ldots, a_{\ell})  = \psi_{\ell-1}^{a_i} (a_1, \ldots, \widehat{a_i}, \ldots, a_{\ell}) \]
where $a_i$ is the rightmost occurrence of the largest entry among $\{a_1, \ldots, a_{\ell}\}$ and the circumflex indicates omission.
\end{theorem}
\begin{proof}
Suppose $\pi(a_1, \ldots, a_{\ell}) = \lambda \in \mathcal{C}_{\ell}^k$.  Then $\lambda$ corresponds to a balanced flush abacus $A$ in which the first row of $\lambda$ corresponds to the last active bead $B$ in the reading order for the abacus, and $B$ occurs on runner $i-1$.  The bijection $\Phi_{\ell}^{k}(\lambda)$ is defined on $A$ by deleting the runner $i-1$ in order to obtain the abacus of an $(\ell-1)$-core.  Since the original abacus is balanced, when we remove runner $i-1$ from $A$ we are left with an abacus $A'$ in which the balance number is $-a_i$.  Applying $\psi_{\ell-1}$ corresponds to shifting all of the entries of $A'$ forward one entry in the reading order of $A'$, or equivalently adding $\mathbf{1}$ to the $\beta$-numbers for $\lambda$.  Hence, applying $\psi_{\ell-1}^{a_i}$ to $A'$ produces a balanced flush abacus for the same partition as $\Phi_{\ell}^{k}(\lambda)$.
\end{proof}

\begin{figure}[p!]  
\centering
\includegraphics{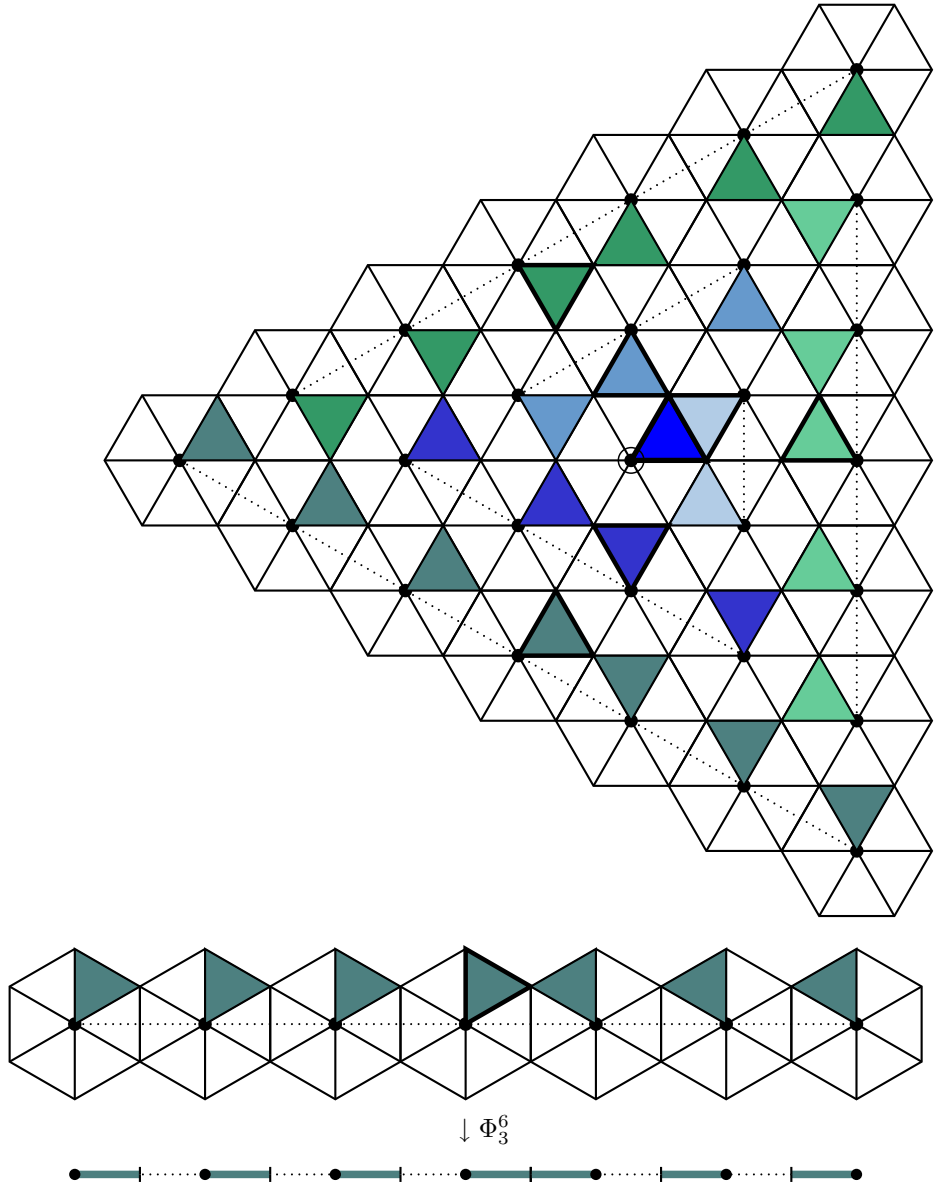}
\caption{\small The bijection $\Phi_{3}^{6}$ as a projection of $A_{2}^{(1)} \cap H_{3}^6 \rightarrow A_{1}^{(1)}$.  The $B_{\circ} + \mathbf{a}$ are hexagons which tile the plane centered at $\mathbf{a}$ which is darkened, and each $B_{\circ} + \mathbf{a}$ contains 6 triangular alcoves.  Each alcove corresponds to an element of $\widetilde{S_{\ell}}$ and each hexagon corresponds to a coset in $\widetilde{S_{\ell}} / S_{\ell}$, so there is a unique 3-core alcove in each hexagon.  These are shaded above.  Two alcoves $B_{\circ} + \mathbf{a}$ and $B_{\circ} + \mathbf{b}$ share the same color only if the first parts of the partitions $\pi(\mathbf{a})$ and $\pi(\mathbf{b})$ agree. }\label{f:sl3}
\end{figure}

The geometric interpretation of $\Phi_{\ell}^k$ as pictured in Figure~\ref{f:sl3} can also be described as follows.  We observed for an $\ell$-core $\lambda$ with $\lambda_1 = k$ that $\pi^{-1}(\lambda)$ lies in the affine hyperplane $H_{\ell}^k$.  We can identify $H_{\ell}^k \cap \Lambda_R$ with a copy of the root lattice of $A_{\ell-2}$ via $(\Phi_{\ell}^k)^{-1}$.  This embedding requires cyclically shifting coordinates as described by $\psi$ and depends on $k$.

\begin{example}\label{e:coords}
Let $\ell = 3$.  The affine hyperplane $H_{3}^7$ contains $(3, 1, -4) = \pi^{-1}(7, 5, 4^2, 3^2, 2^2, 1^2)$.  We decompose $\Phi_{3}^7$ as the composition of a translation and a cyclic shift of coordinates.  Translating by the vector $\mathbf{t} = (-3, 1, 2)$ sends $H_{3}^7$ to the linear hyperplane
\[ V' = \{ (a_1, a_2, a_3) \in V : a_1 = 0 \} \]
and in particular sends $(3, 1, -4)$ to $(0, 2, -2)$.  We view $V'$ as a subspace of $\mathbf{R}^2$ with orthonormal basis $\{e_1', e_2'\}$ and an associated root system of type $A_{\ell-2}$.  Hence, we must identify $e_1'$ with $e_3$ and $e_2'$ with $e_2$.  Therefore, we have $\psi^{3}(1,-4) = (-2, 2)$ corresponding to $\Phi_{3}^7(7,5,4^2,3^2,2^2,1^2) = (4,3,2,1)$.
\end{example}

\begin{remark}
If we focus on $len(\lambda)$ instead of $\lambda_1$, all of the $\ell$-cores with fixed length $m = len(\lambda)$ have $\pi^{-1}(\lambda)$ lying in the affine hyperplane 
\[ \{ \mathbf{a}  \in \mathbf{R}^{\ell} : (\mathbf{a}, \e_{((1-m)\ \mathrm{mod}\ \ell)}) = -\lceil {m \over {\ell} } \rceil \} \cap V. \]
If we drew dotted lines in Figure~\ref{f:sl3} connecting the $\pi^{-1}(\lambda)$ with fixed $m = len(\lambda)$ (instead of those with fixed $k=\lambda_1$), then the lines would appear to spiral backwards from the direction of those in Figure~\ref{f:sl3}.  This can be explained by the fact that sending a partition to its transpose corresponds to the transformation of $\mathbf{R}^{\ell}$ given by $(a_1, \ldots, a_{\ell}) \mapsto (-a_{\ell}, \ldots,- a_1)$ as shown in \cite{garvan-kim-stanton}.
\end{remark}

\subsection{The bijection $\Phi_{\ell}^k$ as a subexpression in Coxeter generators}

Recall from Section~\ref{s:phi_on_diagram} that $\Phi_{\ell}^k$ removes all rows from $\lambda$ in the same equivalence class as the first row.  In Proposition~\ref{p:redexp}, we described a canonical reduced expression $w(\lambda)$ for $\lambda = (\lambda_1, \ldots, \lambda_m) \in \mathcal{C}_{\ell}$.  In this construction, the rows of $\lambda$ are partitioned into $\ell$ equivalence classes which we denote by $[j]$ according to the residue $j$ of their rightmost box.  Let $i$ be the residue of the rightmost box $\mathsf{B}$ in the last row $\lambda_m$ of $\lambda$.  Observe that $\mathsf{B}$ is a removable $i$-box unless $\lambda = \emptyset$, and applying $s_i$ removes one box from each of the rows $\equiv [i]$.

We claim that two rows of $s_i (\lambda)$ are equivalent if and only if the rows were equivalent in $\lambda$.  To see this, consider that there can be no rows $\equiv [i-1]$ in $\lambda$.  Otherwise there exists a box in the same column as $\mathsf{B}$ whose row is $\equiv [i-1]$, and so the hooklength of this box is divisible by $\ell$ which contradicts $\lambda$ being an $\ell$-core.  Therefore, the rows $\equiv [i-1]$ in $s_i(\lambda)$ are precisely the rows $\equiv [i]$ in $\lambda$ with the possible exception of $\lambda_m$ if $\lambda_m = 1$.  In any case, no other rows change equivalence classes.

Suppose $w(\lambda) = s_{i_1} s_{i_2} \cdots s_{i_p}$ is the canonical reduced expression for $\lambda$ obtained from Proposition~\ref{p:redexp}, so $i_1$ is the residue of the box $\mathsf{B}$ and $i_p = 0$.  Working left to right to reduce $w(\lambda)$ to the identity, each application of $s_{i_j}$ removes a box from every row that is in the same equivalence class as the last row of the intermediate partition.  Let $J$ be the subset of $\{1, \ldots, p\}$ such that the first row and the last row are not in the same equivalence class in the $\ell$-core $s_{i_j} s_{i_{j+1}} \cdots s_{i_p}(\emptyset)$.  Then the subexpression of $w(\lambda)$ corresponding to the indices in $J$ gives the canonical minimal length coset representative for the $(\ell-1)$-core $\Phi_{\ell}^k(\lambda)$ after relabeling the residues with respect to $\ell-1$.

Since we remove a box from every row in the equivalence class of the last row at each step, we remove in particular a box from the longest row in that equivalence class.  Hence, we see that the positions $\{1, \ldots, p\} \setminus J$ that are deleted from $w(\lambda)$ in the application of $\Phi_{\ell}^k$ correspond with boxes in the first row of $\lambda$, so applying $\Phi_{\ell}^k$ reduces the Coxeter length by exactly $k$.

\begin{example}
Suppose $\ell = 5$ and let $\lambda=(9,5,3,2,2,1,1,1,1)$.  
We label the diagram of $\lambda$ as shown in Figure~\ref{f:subexp} by residues with respect to $\ell$ on the bottom of each box, and with respect to $\ell-1$ on the top of those boxes that are not in rows equivalent to the first row.  We put dots as placeholders in the tops of these boxes.

\begin{figure}[h]
\[
\lambda = \tableau{ { }^{.}_{0} &  { }^{.}_{1} &  { }^{.}_{2} &  { }^{.}_{3} &  { }^{.}_{4} &  { }^{.}_{0} &  { }^{.}_{1} &  { }^{.}_{2} &  { }^{.}_{3} \\
 { }^{.}_{4} &  { }^{.}_{0} &  { }^{.}_{1} &  { }^{.}_{2} &  { }^{.}_{3} \\
 { }^{0}_{3} &  { }^{1}_{4} &  { }^{2}_{0} \\
 { }^{.}_{2} &  { }^{.}_{3} \\
 { }^{3}_{1} &  { }^{0}_{2} \\
 { }^{2}_{0} \\
 { }^{1}_{4} \\
 { }^{.}_{3} \\
 { }^{0}_{2} }
\]
\caption{$\Phi_{\ell}^k$ as a Coxeter subexpression of $w = {\bf s_2} s_3 {\bf s_4 s_0 s_1} s_2 {\bf s_4 s_3 } s_1 s_0 s_4 s_3 s_2 s_1 s_0$.} \label{f:subexp}
\end{figure}

We see directly from this diagram that 
\[ w(\lambda) = w = {\bf s_2} s_3 {\bf s_4 s_0 s_1} s_2 {\bf s_4 s_3 } s_1 s_0 s_4 s_3 s_2 s_1 s_0 \]
where the entries in $J$ have been highlighted.  For example, applying $s_2$ to $\lambda$ removes the last row as well as the last box from row 5.  Since we did not remove a box from the first row, we have position $1 \in J$.  The next step applies $s_3$ to remove the last row of $s_2(\lambda)$ as well as the last box from rows 1, 2, and 4, so position $2 \notin J$.

After shifting the residues in the bold subexpression $s_2 s_4 s_0 s_1 s_4 s_3$ to be calculated relative to $\ell-1$ as shown in the diagram, we obtain $w(\Phi_{\ell}^{k}(\lambda)) = s_0 s_1 s_2 s_3 s_1 s_0$.  The Coxeter length has been reduced by $15 - 6 = 9 = \lambda_1$.
\end{example}

\section{Relation with a correspondence of Lapointe--Morse}\label{s:l-m}

\subsection{$\Phi_{\ell}^k$ interpreted in terms of the Lapointe--Morse correspondence}
In this section, we show that the bijection $\Phi_{\ell}^k$ can be succinctly expressed as a map between $k$-bounded partitions (one whose first part is $\leq k$) and $(k-1)$-bounded partitions, using the correspondence
\[ \rho_{k+1} : \{ (k+1)\text{-cores} \} \rightarrow \{ \text{partitions with first part } \leq k \} \] 
of Lapointe and Morse \cite[Section 3]{LM} to which we refer the reader for more details.  Let $\widetilde{\Phi_{\ell}^k}$ be defined by the property that the left square in the following diagram commutes, and define $\Upsilon_{\ell}^k$ to simply delete the first column from the diagram of the partition.  Then, we claim that the right square in the diagram also commutes.  Here, $tr$ denotes the transpose of a partition.

\medskip

\[
\begin{CD}
\left\{ \parbox{0.45in}{$\lambda \in \mathcal{C}_{\ell}^k$} \right\} @> tr >> \left\{ \parbox{1.1in}{$\lambda \in \mathcal{C}_{\ell}$, $len(\lambda) = k$} \right\} @> \rho_{\ell} >> \left\{ \parbox{1in}{partitions $\nu$ with $\nu_1 \leq \ell-1$ and $len(\nu) = k$} \right\} \\
@ V{\Phi_{\ell}^k} VV @ VV {\widetilde{\Phi_{\ell}^k}} V @ VV {\Upsilon_{\ell}^k} V \\
\left\{ \parbox{0.55in}{$\mu \in \mathcal{C}_{\ell-1}^{\leq k}$} \right\} @> tr >> \left\{ \parbox{1.25in}{$\mu \in \mathcal{C}_{\ell-1}$, $len(\mu) \leq k$} \right\} @>> \rho_{\ell-1} > \left\{ \parbox{1in}{partitions $\sigma$ with $\sigma_1 \leq \ell-2$ and $len(\sigma) \leq k$} \right\} \\
\end{CD}
\]

\bigskip

Recall that $h_{\mathsf{B}}^{\lambda}$ denotes the hooklength of a box $\mathsf{B}$ in a partition diagram $\lambda$.  Let $\lambda = (\lambda_1, \ldots, \lambda_k)$ be an $\ell$-core with first column length $k$.  Observe that $\widetilde{\Phi_{\ell}^k}$ acts by removing entire \em columns \em from $\lambda$ by transposing the explanation in Section~\ref{s:phi_on_diagram}.  Row-wise, we can describe $\widetilde{\Phi_{\ell}^k}$ as removing the leftmost box $\mathsf{B}$ in the row together with all of the boxes $\mathsf{B}'$ in the row having $h_{\mathsf{B}}^{\lambda} \equiv h_{\mathsf{B}'}^{\lambda} \mod \ell$.

As in \cite{LM}, we view $\rho_{\ell}(\lambda)$ as the result of left-justifying all of the rows in a skew-diagram $\lambda / \gamma$, where $\gamma$ consists of the boxes of $\lambda$ having hooklength $> \ell$.  We say that the boxes of $\lambda$ lying in $\lambda / \gamma$ are the \em skew boxes \em while the boxes of $\gamma$ are the \em non-skew boxes\em.  

\begin{example}
Suppose $\ell = 5$ and consider the $\ell$-core $\lambda = (6,4,3,3,2,1,1)$.  In the diagrams below, we have labeled the boxes by their hooklengths.  The skew-boxes of $\gamma \subset \lambda$ are indicated in boldface, so $\rho_{\ell}(\lambda) = (3,2,2,2,2,1,1)$.  The entries that are deleted in the application of $\widetilde{\Phi_{\ell}^k}$ are indicated with underline.  
\[
\lambda = \tableau{ \underline{12} & 9 & \underline{7} & {\bf 4} & \underline{\bf 2} & {\bf 1} \\
\underline{9} & 6 & \underline{\bf 4} & {\bf 1} \\
\underline{7} & {\bf 4} & \underline{\bf 2} \\
\underline{6} & {\bf 3} & \underline{\bf 1} \\
\underline{\bf 4} & {\bf 1} \\
\underline{\bf 2} \\
\underline{\bf 1} \\
}
\ \ \hspace{0.4in} \ \ \
\widetilde{\Phi_{\ell}^k}(\lambda) = \tableau{ 7 & {\bf 3} & {\bf 1} \\
5 & {\bf 1} \\
{\bf 3} \\
{\bf 2} \\
{\bf 1} \\
}
\]
From these diagrams, we see that $\widetilde{\Phi_{\ell}^k}(\lambda)$ is a 4-core and $\Upsilon_{\ell}^k(\rho_{\ell}(\lambda)) = (2,1,1,1,1) = \rho_{\ell-1}(\widetilde{\Phi_{\ell}^k}(\lambda))$.
\end{example}

To simplify notation, let $\widetilde{\lambda} = \widetilde{\Phi_{\ell}^k}(\lambda)$ and define $\widetilde{\gamma}$ to consist of the boxes of $\widetilde{\lambda}$ having hooklength $> \ell-1$ so that $\rho_{\ell-1}(\widetilde{\lambda})$ is the result of left justifying the boxes of $\widetilde{\lambda} / \widetilde{\gamma}$.  

\begin{lemma}\label{l:lm-one-rem}
There is exactly one skew box deleted from each row of $\lambda$ in the application of $\widetilde{\Phi_{\ell}^k}$.  Also, a box $\mathsf{B}$ in $\lambda$ that is not deleted in the application of $\widetilde{\Phi_{\ell}^k}$ is skew with respect to $\ell$ if and only if the corresponding box $\widetilde{\mathsf{B}}$ of $\widetilde{\Phi_{\ell}^k}(\lambda)$ is a skew box with respect to $\ell-1$.
\end{lemma}

\begin{proof}
It suffices to prove these statements for a fixed row.  Since the skew boxes are those with hooklength $< \ell$, we delete at most one skew box from each row of $\lambda$ when we apply $\widetilde{\Phi_{\ell}^k}$.  Next, we show that at least one skew box is deleted.  Let $\mathsf{L}$ be the leftmost box of $\lambda$ in row $i$.  Since $\lambda$ is an $\ell$-core, the partition $\widehat{\lambda} = (\lambda_i, \lambda_{i+1}, \ldots, \lambda_k)$ is also an $\ell$-core.  Form an unbalanced abacus $A$ of $\widehat{\lambda}^{tr}$ such that the beads of $A$ correspond to hooklengths of the first \em row \em of $\widehat{\lambda}$ which is the $i^{\textrm{th}}$ row of $\lambda$.  In particular, those beads corresponding to the boxes with equivalent hooklength to $h_{\mathsf{L}}^{\lambda}$ form the rightmost longest runner of the abacus $A$ and they are flush.  Hence, there exists a box $\mathsf{L}'$ in the same row as $\mathsf{L}$ having hooklength $\equiv h_{\mathsf{L}}^{\widehat{\lambda}} = h_{\mathsf{L}}^{\lambda} \mod \ell$ such that $1 \leq h_{\mathsf{L}'}^{\lambda} \leq \ell-1$, so the corresponding box $\mathsf{L}'$ is a skew box.  The box $\mathsf{L}'$ is deleted from $\lambda$ when we apply $\widetilde{\Phi_{\ell}^k}$.

To prove the second statement, suppose $\mathsf{B}$ is a box in the $i$th row of $\lambda$ that does not get deleted in $\widetilde{\lambda}$.  Let $\widetilde{\mathsf{B}}$ be the corresponding box in $\widetilde{\lambda}$.  Then $\mathsf{B}$ corresponds to an active bead in the abacus $A$ that is not on same runner as $\mathsf{L}$.  We have that $\mathsf{B}$ is skew with respect to $\ell$ if and only if the corresponding bead lies on level 0 of the abacus $A$.  Since removing the runner containing $\mathsf{L}$ does not change any levels of the remaining beads, we have that $\widetilde{\mathsf{B}}$ is skew with respect to $\ell-1$ if and only if $\mathsf{B}$ is skew with respect to $\ell$.
\end{proof}

\begin{theorem}
The map $\Upsilon_{\ell}^k$ is a bijection which makes the diagram at the beginning of this section commute.
\end{theorem}
\begin{proof}
We see from Lemma~\ref{l:lm-one-rem} that the notion of skew box is preserved under the application of $\widetilde{\Phi_{\ell}^k}$.  Thus, we have that $\rho_{\ell-1}(\widetilde{\Phi_{\ell}^k}(\lambda))$ is formed from $\rho_{\ell}(\lambda)$ by simply deleting the first column.
\end{proof}

\bibliographystyle{amsalpha}

\end{document}